\documentclass[leqno,10pt]{article}

\input{amssym}
\begin{document}
\title{ Symmetry Analysis for a New Form of
the\\ Vortex Mode Equation}
\date{}
\author{Mehdi Nadjafikhah \and Ali Mahdipour--Shirayeh}
\maketitle
\begin{abstract} Giving a new form of the vortex mode equation by a proper change
of parameter, our aim is to analyze the point and contact
symmetries of the new equation. Fundamental invariants and a form
of general solutions of point transformations along with some
specific examples are also derived.
\end{abstract}
\noindent {\bf M.S.C. 2000:} 34C14, 58D19, 35L05.
\section{Introduction}
Investigation of nonlinear phenomena appearing in a very wide area
of pure and applied sciences has met extremely extensive
progresses and developments. These studies which split into
numerical and analytical considerations are essentially and in
most cases related to some nonlinear differential equations. Among
those nonlinear systems, a few interesting open problems concern
the hydrodynamic type of equations governing fluid motions.
Especially the Euler and Navier-Stokes equations which reveal a
mysterious behavior are being intensively studied in two main
considerations: The incompressible motion mostly dealing with
vortex dynamics and the compressible flow concerning the
appearance of discontinuities shocks (see \cite{Esh-Ma} and
related references therein). In \cite{Esh-Ma} after a brief
derivation of relativistic ideal fluid equations, a
multidimensional simple wave ansatz is substituted into these
equations and various modes (for instance the vortex mode) and
phase velocities relative to the laboratory (fixed) frame are
found.

The vortex mode equation \cite{Boil,Esh-Ma,Web-Rat} is defined as
a first order ODE
\begin{eqnarray}
\frac{d{\bf k}}{d\varphi}\cdot\Big({\bf n}-\frac{{\bf k}\cdot{\bf
n}}{{\bf k}^2+w}\,{\bf k}\Big)=0,\label{eq:1}
\end{eqnarray}
where $w$ is a constant, $\varphi$ is still treated as the wave
phase, and ${\bf k}=(k_1,k_2,k_3)$ and ${\bf n}=(n_1,n_2,n_3)$ are
some vectors in ${\Bbb R}^3$ with the physical meaning of
${\bf\kappa}$ in \cite{Esh-Ma} and unit normal vector to the wave
front resp. A symmetry analysis of Eq.~(\ref{eq:1}) up to both
point and contact transformations has performed in \cite{Esh-Ma}.
But in this paper, we investigate symmetry properties of a new
form of Eq.~(\ref{eq:1}) as a simple form of the vortex mode
equation. As is well known, under change of coordinates the
symmetry group of a system of differential equations remains
unchanged. But since the jacobian of the following change of
parameter is zero, so symmetry analysis of the vortex mode
equation and the new form are not necessarily the same. By
applying the following change of parameter
\begin{eqnarray*}
t:=\frac{1}{2}\,\ln({\bf k}^2+w),
\end{eqnarray*}
we find the new form as follows
\begin{eqnarray}
{\bf n}\cdot\Big(\frac{d{\bf k}}{dt}-{\bf k}\Big)=0.\label{eq:2}
\end{eqnarray}
Since Eq.~(\ref{eq:1}) is a homogeneous linear equation with
respect to ${\bf n}$, so we consider it to be of arbitrary length
and not necessarily unit.

Eq.~(\ref{eq:2}) is in fact an expression of the vortex mode
equation that provides an in depth study of Eq.~(\ref{eq:1}).
Roughly speaking, it leads to slightly simpler calculations for
finding exact solutions of the vortex mode equation. But for
reaching to this goal, we investigate symmetry properties of
Eq.~(\ref{eq:2}) which play a key role in finding general
solutions, fundamental invariants, invariant solutions and etc.
Moreover, knowledge of a symmetry group of Eq.~(\ref{eq:2}) allows
us to construct new solutions from old ones \cite{Ib,Ib2,Ol,Ol2}.
Therefore in this study, we concern with the latter equation to
find its point and contact symmetry properties and also give its
fundamental invariants and a form of general solutions.

Throughout this paper we assume that indices $i,j$ varies between
1 and 3 and each index of a function implies the derivation of the
function with respect to it, unless specially stated otherwise.
%%%%%%%%%%%%%%%%%%%%%%%%%%%%%%%%%%%%%%%%%%%%
\section{The point Symmetry of the Equation}

To find the symmetry group of Eq.~(\ref{eq:2}) by Lie
infinitesimal method, we follow the method presented in \cite{Ol}.
We find infinitesimal generators of the equation and also the Lie
algebra structure of the symmetry group of (\ref{eq:2}). In this
section, we are concerned with the action of the point
transformation group.

The equation is a relation among with the variables of 1--jet
space $J^1({\Bbb R},{\Bbb R}^6)$ with (local) coordinate $(t, {\bf
k}, {\bf n}, {\bf q}, {\bf p})=(t, k_i, n_j,  q_r, p_s)$
($\mbox{for}\;\;1\leq i,j,r,s\leq 3$), where this coordinate
involving a independent variable $t$ and 6 dependent variables
$k_i,n_j$ and their derivatives $q_r,p_s$ of first order with
respect to $t$ resp.

Let ${\cal M}$ be the total space of independent and dependent
variables. The solution space of Eq.~(\ref{eq:2}), (if it exists)
is a subvariety $S_{\Delta}\subset J^1({\Bbb R},{\Bbb R}^6)$ of
the first order jet bundle of one--dimensional submanifolds of
${\cal M}$.

We define a point transformation on ${\cal M}$ with relations
\begin{eqnarray*}
\tilde{t}=\phi(t,k_i,n_j),\hspace{1cm}
\tilde{k}_r=\chi_r(t,k_i,n_j),\hspace{1cm}
\tilde{n}_s=\psi_s(t,k_i,n_j).
\end{eqnarray*}
for $\phi,\chi_r$ and $\psi_s$ are some smooth functions. Let
$$v:=T\,\displaystyle{\frac{\partial }{\partial
t}}+\sum_{i=1}^3\Big( K_i\displaystyle{\frac{\partial }{\partial
k_i}}+ N_i\displaystyle{\frac{\partial }{\partial n_i}}\Big)$$ be
the general form of infinitesimal generators that signify the Lie
algebra ${\frak g}$ of the symmetry group $G$ of Eq.~(\ref{eq:2}).
In this relation, $T,K_i$ and $N_j$ are smooth functions of
variables $t,k_i$ and $n_j$. The first order prolongation
\cite{Ol} of $v$  is as follows
\begin{eqnarray}
v^{(1)}&:=& v + \sum_i\,K_i^t\frac{\partial}{\partial q_i}+
\sum_j\,N_j^t\frac{\partial}{\partial p_j} ,\nonumber
\end{eqnarray}
where $K_i^t=D_t\,Q_1^i+T\,q_{i,t}$ and
$N_j^t=D_t\,Q_2^j+T\,p_{j,t}$, in which $D_t$ is total derivative
and $Q_1^i= K_i - T\,q_i$ and $Q_2^j= N_j - T\,p_j$ are
characteristics of vector field $v$ \cite{Ol}. By effecting
$v^{(1)}$ on (\ref{eq:2}), we obtain the following expression
\begin{eqnarray}
&& \sum_i \Big[ n_i\,(K_{i t} - K_i)  - N_i\,k_i + p_i\,\sum_j
n_j\,K_{jn_i} + q_i\Big( N_i - n_i\,T_t + \sum_j\,
n_j\,K_{jk_i}\Big)\nonumber\\[-3.5mm]
&&  \label{eq:3}\\[-2.5mm]
&& - q_i^2\,n_i\,T_{k_i}\Big]  - \sum_{i\neq j}
q_i\,q_j(n_i\,T_{k_j}+ n_j\,T_{k_i}) - n_i\,\sum_{i,j}
q_i\,p_j\,T_{n_j} = 0,\nonumber
\end{eqnarray}
in which, each index (exception for determined indices) signifies
the derivation with respect it.\par
We can prescribe $t, k_i, n_j, q_r, p_ s$ ($1\leq i,j,r,s\leq3$)
arbitrarily, and functions $K_i$ and $N_j$ only depend on
$t,k_i,n_j$ ($i,j=1,2,3$). So, Eq.~(\ref{eq:3}) will be satisfied
if and only if we have the following equations
\begin{eqnarray}
&& N_i - n_i\,T_t + n_1K_{1k_i}+ n_2K_{2k_i} + n_3K_{3k_i}= 0, \label{eq:4}\\
&& n_1K_{1n_i}+n_2K_{2n_i}+n_3K_{3n_i}= 0, \label{eq:5}\\
&& n_i\,T_{k_i}=0, \hspace{0.5cm}n_i\, T_{n_j}= 0, \hspace{0.5cm}n_i\,T_{k_j}+n_j\,T_{k_i}=0, \label{eq:6}\\
&&  \sum_{i=1}^3
\Big(n_i(K_{i\,t}-K_i)-N_i\,k_i\Big)=0,\label{eq:7}
\end{eqnarray}
when $1\leq i,j\leq 3$ and (\ref{eq:2}) is satisfied. These
equations are called the determining equations. From (\ref{eq:5})
for each $i$, we have
\begin{eqnarray}
N_i = n_i\,T_t - (n_1K_1+ n_2K_2 + n_3K_3)_{k_i}.\label{eq:8}
\end{eqnarray}
Since ${\bf n}\neq0$, without less of generality, one may assume
that $n_1\neq0$. Then by Eqs. (\ref{eq:6}) we conclude that $T$
just depends on $t$:
\begin{eqnarray*}
T= T(t).
\end{eqnarray*}
By solving the Eqs. (\ref{eq:5}) along with Eq.~(\ref{eq:7}) in
respect to $K_1$, $K_2$ and $K_3$, then we deduce the following
relations (provided by \textsc{Maple})
\begin{eqnarray}
&& \hspace{-1cm}K_1 = k_1\,T + e^t\,F^1  \label{eq:9}\\
&&\hspace{-1cm} \nonumber\\
&&\hspace{-1cm} K_{2t} = K_2-k_1\,K_{2k_1}-k_2\,K_{2k_2}-
k_3\,K_{2k_3}+ k_2 \,T_t,  \nonumber\\[-3mm]
&&\hspace{-1cm}  \label{eq:10} \\[-2mm]
&&\hspace{-1cm} n_1\,K_{2n_1} =e^t\,F^1_4, \hspace{0.5cm}
n_2\,K_{2n_2} = -e^t\,F^1_4 -n_3\,K_{2n_3},
\nonumber\\[1mm]
&&\hspace{-1cm}\nonumber\\
&&\hspace{-1cm} K_{3t} = K_3-k_1\,K_{3k_1}-k_2\,K_{3k_2}- k_3\,K_{3k_3}+ k_3 \,T_t, \nonumber\\[-3mm]
&&\hspace{-1cm} \label{eq:11} \\[-3mm]
&&\hspace{-1cm} n_1\,K_{3n_1} = e^t\,F^1_5,\hspace{0.5cm} K_{3n_2}
= K_{2n_3}\hspace{0.5cm} n_3\,K_{3n_3} = -e^t\,F^1_5 -
n_2\,K_{2n_3}\nonumber
\end{eqnarray}
for arbitrary function
$F^1=F^1(k_1e^{-t},k_2e^{-t},k_3e^{-t},\frac{n_2}{n_1},\frac{n_3}{n_1})$.
In these relations, $F^1_i$ implies the derivation of $F^1$ in
respect to the $i^{th}$ coefficient.\par
Eqs. (\ref{eq:10}) lead to the below relations
\begin{eqnarray*}
&& K_{2t}+ k_1\,K_{2k_1}+ k_2\,K_{2k_2}+ k_3\,K_{2k_3}-K_2 -
k_2\,T_t=0,\\
&& n_1\,K_{2n_1}+ n_2\,K_{2n_2}+ n_3\,K_{2n_3}=0,
\end{eqnarray*}
so that after solving determines the form of $K_2$ as
\begin{eqnarray}
K_2= k_2\,T + e^t\,F^2,\label{eq:12}
\end{eqnarray}
when
$F^2=F^2(k_1e^{-t},k_2e^{-t},k_3e^{-t},\frac{n_2}{n_1},\frac{n_3}{n_1})$
is an arbitrary smooth function. Also, from Eqs. (\ref{eq:11}), we
find that
\begin{eqnarray*}
&& K_{3t}+ k_1\,K_{3k_1}+ k_2\,K_{3k_2}+ k_3\,K_{3k_3}- K_3 -
k_3\,T_t=0,\\
&& n_1\,K_{3n_1}+ n_2\,K_{3n_2}+ n_3\,K_{3n_3}=0,
\end{eqnarray*}
which these expressions tend to the following solution of $K_3$
with respect to arbitrary smooth function
$F^3=F^3(k_1e^{-t},k_2e^{-t},k_3e^{-t},\frac{n_2}{n_1},\frac{n_3}{n_1})$:
\begin{eqnarray}
K_3= k_3\,T + e^t\,F^3.\label{eq:13}
\end{eqnarray}
But by satisfying $F^2$ and $F^3$ resp. in the two last relations
of (\ref{eq:10}) and three last relations of (\ref{eq:11}), we
find that
\begin{eqnarray*}
&&\hspace{-0.7cm} F^1 =
-\frac{1}{n_1}(n_2\,F^2+n_3\,G)-n_2\!\int\!\frac{1}{n_1^2}F^2\,dn_1
- n_3\!\int\!\frac{1}{n_1^2}G\,dn_1+
\frac{n_2n_3}{n_1}\!\int\!\frac{1}{n_1^2}F^2_5\,dn_1 \\
&& +\frac{1}{n_1}\!\int\!\Big(F^2-\frac{n_3}{n_1}\,F^2_5\Big)\,dn_2 + \frac{1}{n_1}\int G\,dn_3+ H,\\
&&\hspace{-0.7cm} F^3 = \frac{1}{n_1}\int\!F^2_5\, dn_2 -
n_2\int\!\frac{1}{n_1^2}F^2_5\, dn_1 + G,
\end{eqnarray*}
when $F^2$, $G=G(k_1e^{-t},k_2e^{-t},k_3e^{-t},\frac{n_3}{n_1})$
and $H=H(k_1e^{-t},k_2e^{-t},k_3e^{-t})$ are arbitrary functions
and $F^j_i$ denotes the derivation of $F^j$ with respect to its
$i^{th}$ coefficient (similar statement is valid for $G_i$ and
$H_i$ in subsequent relations).\par
Substituting the new forms of $K_1, K_2$ and $K_3$ in Eqs.
(\ref{eq:5}), we obtain the following relations
\begin{eqnarray}
&&\hspace{1cm} \frac{n_2n_3}{n_1}\int\frac{1}{n_1^2}F^2_5\,dn_1+\frac{1}{n_1}(n_2\,F^2+n_3\,G)=0, \label{eq:14}\\
&&\hspace{1cm} F^2 + n_1n_3\int\frac{1}{n_1^3}\,F^2_5\,dn_1=0, \label{eq:15}\\
&&\hspace{1cm} n_2(n_1+n_3)F^2_5 +n_2^2\,F^2_4  +
n_1^3n_2\int\frac{1}{n_1^3}F^2_5\,dn_1 - n_1^2\,G = n_1^2n_2\int
\frac{1}{n_1^2}F^2_5\,dn_1. \label{eq:16}
\end{eqnarray}
By applying the latest relations in $F^1$ and $F^3$, we attain the
below relations (for arbitrary $F^2$)
\begin{eqnarray}
&& F^1 =
-\frac{2}{n_1}(n_2\,F^2+n_3\,G)-n_2\!\int\!\!\frac{1}{n_1^2}F^2\,dn_1+\frac{1}{n_1}
\!\!\int\!\Big(F^2-\frac{n_3}{n_1}\,F^2_5\Big)\,dn_2 \nonumber\\
&&\hspace{1cm} -n_3\!\int\!\!\frac{1}{n_1^2}G\,dn_1+ \frac{1}{n_1}\int G\,dn_3+ H,  \label{eq:17}\\
&& F^3 = \frac{1}{n_1}\int\!F^2_5\, dn_2 + \frac{n_2}{n_3}F^2 +
2\,G,\label{eq:18}
\end{eqnarray}
in which, we assumed that $n_3\neq 0$. Otherwise, from
Eq.~(\ref{eq:15}), $F^2=0$ and hence by Eq.~(\ref{eq:16}), $G=0$
and therefore we have $F^1=H$ and $F^3=0$.\par
We continue our investigation of symmetry group with the condition
$n_3\neq 0$. From Eq.~(\ref{eq:14})--(\ref{eq:16}), we can elicit
the following relation
\begin{eqnarray}
n_2(n_1+n_3)\,F^2_5 + n_2^2\,F^2_4=0.\label{eq:19}
\end{eqnarray}
On the other hand, by replacing $F^2_5=F^3_4$ from relation
$K_{3n_2}=K_{2n_3}$ of (\ref{eq:11}), in (\ref{eq:19}), we find
that $n_2\,F^2=-2\,n_3\,G$. This relation along with relation
(\ref{eq:19}) leads to an equation that its solution determines
the form of $G$ as following
\begin{eqnarray*}
G=\frac{1}{n_3}(n_1+n_3)\,L,
\end{eqnarray*}
where $L=L(k_1e^{-t},k_2e^{-t},k_3e^{-t})$ is an arbitrary smooth
function. Hence, $F^2=-\frac{2}{n_2}(n_1+n_3)\,L$ and from Eqs.
(\ref{eq:17}) and (\ref{eq:18}), we have (we suppose that $n_2\neq
0$, otherwise from Eq.~(\ref{eq:14}), $G=0$ and hence $L=0$)
\begin{eqnarray*}
F^1 = \frac{1}{n_1}\Big(2\,(n_1+n_2+n_3) +
n_1\ln(n_1\,n_3)\Big)\,L +H,\hspace{1.65cm} F^3 = -2\,\ln(n_2)\,L.
\end{eqnarray*}
Applying the last forms of $K_i\,$s in (\ref{eq:5}) for $i=3$, we
find $\frac{n_1}{n_3}\,L=0$. We have assumed that $n_1\neq 0$,
therefore $L=0$ and the forms of $K_i$~s and $N_j$~s
\begin{eqnarray*}
&&K_1 = k_1\,T+ e^t\,H,\hspace{1cm} K_2 = k_2\,T,\hspace{1cm}  \\
&&K_3 = k_1\,T,\hspace{2.1cm}  N_j = n_i(T_t-T) - n_1\,H_i,
\end{eqnarray*}

and finally, the general form of infinitesimal generators as
elements of point symmetry algebra of Eq.~(\ref{eq:2}), which we
call them as {\it point infinitesimal generators}, is as the
following relation, for arbitrary functions $T$ and $H$
\begin{eqnarray}
v &=& T\,\Big(\frac{\partial}{\partial
t}+k_2\,\frac{\partial}{\partial k_2} +
k_3\,\frac{\partial}{\partial k_3}\Big) + (k_1\,T + e^t\,H
)\frac{\partial}{\partial k_1}  \nonumber\\
&&+ \sum_{i=1}^3\Big(n_i(T_t-T) -
n_1\,H_i\Big)\frac{\partial}{\partial n_i}.\label{eq:20}
\end{eqnarray}
One may divides $v$ into the following infinitesimal generators
\begin{eqnarray}
&&\hspace{-0.7cm} v_T = T\Big(\frac{\partial}{\partial t}+
\sum_{i=1}^3 k_i\,\frac{\partial}{\partial k_i}\Big)
+(T_t-T)\sum_{i=1}^3n_i\frac{\partial}{\partial
n_i}, \nonumber\\[-3.5mm]
&&\label{eq:21}\\[-3.5mm]
&&\hspace{-0.7cm}  v_H = e^t\,H\,\frac{\partial}{\partial k_1}-
n_1\,\sum_{i=1}^3\,H_i\frac{\partial}{\partial n_i}.\nonumber
\end{eqnarray}
The Lie bracket (commutator) of vector fields (\ref{eq:21})
straightforwardly a linear combination of them. The table of
commutators is given in Table 1. Hence, the Lie algebra ${\frak
g}=\langle v_T,v_H \rangle$ of point symmetry group $G$ is an
abelian Lie algebra.
\begin{table}[t]
\begin{eqnarray*}\begin{array}{l|l l l l l l l} \hline\hline
[\,,\,]  &\hspace{1cm}v_T\hspace{0.5cm}& v_H\\ \hline
v_T      &\hspace{1cm} 0          \hspace{2cm}&  0 \hspace{1cm}     \\[1mm]
v_H      &\hspace{1cm} 0          \hspace{2cm}&  0 \hspace{1cm}     \\[1mm]
\hline\hline
\end{array}\end{eqnarray*}
\noindent\caption{The commutators table of ${\frak g}$ for
Eq.~(\ref{eq:2})}\label{table:1}
\end{table}
\paragraph{Theorem 1.}{\it The set of all point infinitesimal generators in the forms of
(\ref{eq:21}) is the infinite dimensional abelian Lie algebra of the point symmetry group
of equation (\ref{eq:2}).}\\

According to theorem 2.74 of \cite{Ol}, the invariants
$u=I(t,k_1,k_2,k_3,n_1,n_2,n_3)$ of one--parameter group with
infinitesimal generators in the form of (\ref{eq:21}) satisfy the
linear, homogeneous partial differential equations of first order:
\begin{eqnarray*}
v[I]=0.
\end{eqnarray*}
The solutions of the latter, are found by the method of
characteristics (See \cite{Ol} and \cite{Ib} for details). So we
can replace the latest equation by the following characteristic
system of ordinary differential equations ($1\leq i,j\leq 3$)
\begin{eqnarray}
\frac{dt}{T}=\frac{dk_i}{K_i}=\frac{dn_j}{N_j}.\label{eq:22}
\end{eqnarray}
By solving the Eqs. (\ref{eq:22}) of the differential generator
(\ref{eq:21}), we (locally) find the following general solutions
\begin{eqnarray}
&&\hspace{0.5cm} I_1(t,{\bf k},{\bf n}) = k_2^{-1}(k_1\,T + H)=d_1, \nonumber\\
&&\hspace{0.5cm} I_i(t,{\bf k},{\bf n}) = \ln(k_i)-t
=d_i, \hspace{0.25cm}(\mbox{for} \hspace{0.2cm} i=2,3) \nonumber\\[-2.5mm]
&&\hspace{0.5cm} \label{eq:23}\\[-2.5mm]
&&\hspace{0.5cm} I_4(t,{\bf k},{\bf n}) = (T_t-T-H_1)\,\ln(k_1)-T\,\ln(n_1) = d_4, \nonumber\\
&&\hspace{0.5cm} I_j(t,{\bf k},{\bf n}) =
(T_t-T)\,\ln(k_j)-T\,\ln\left(n_j(T_t-T)-n_1\,H \right)= d_j,
\hspace{0.25cm}(\mbox{for} \hspace{0.2cm} j=5,6). \nonumber
\nonumber
\end{eqnarray}
when $d_i$~s are some constants. The functions
$I_1,I_2,\cdots,I_6$ form a complete set of functionally
independent invariants of one--parameter group generated by
(\ref{eq:21}) (see \cite{Ol}).\par
Similar to the theorem of section 4.3.3 of \cite{Ib}, the derived
invariants (\ref{eq:18}) as independent first integrals of the
characteristic system of the infinitesimal generator
(\ref{eq:16}), provide the general solution $$S(t,{\bf k},{\bf
n}):=\mu(I_1(t,{\bf k},{\bf n}),I_2(t,{\bf k},{\bf
n}),\cdots,I_6(t,{\bf k},{\bf n})),$$ with an arbitrary function
$\mu$, which satisfies in the equation $v[\mu]=0$.

This theorem can be extended for each finite set of independent
first integrals (invariants) of characteristic system provided
with an infinitesimal generator.\par
In the following, we give some examples provided with different
selections of coefficients of Eq.~(\ref{eq:20}) for better
studying, and we assume that each appeared coefficient of vector
fields be non-zero.
\paragraph{Example 1.} If we assume that $T=1$ and $H=0$, then the infinitesimal
operator (\ref{eq:15}) reduces to the following vector field
\begin{eqnarray*}
v_1=\frac{\partial}{\partial t}+\sum_{i=1}^3
k_i\,\frac{\partial}{\partial k_i} -\sum_{j=1}^3
n_j\,\frac{\partial}{\partial n_j}.
\end{eqnarray*}
and the group transformations (or flows) for the parameter $s$ are
expressible as $(t,k_i,n_j)\rightarrow
(t+s,k_i\,e^s,n_j\,e^{-s})$, that form the (local) symmetry group
of $v_1$.\par
The derived invariants in this case will be as follows
\begin{eqnarray*}
I_i=\ln(k_i)-t,\hspace{1cm} I_{j+3}=\ln(n_j)+t,\hspace{1cm}
i,j=1,2,3.
\end{eqnarray*}
Therefore, the general solution corresponding to $v_1$, when $\mu$
is an arbitrary function, will be $S(t,{\bf k},{\bf n})=
\mu\Big(\ln(k_i)-t,\ln(n_j)+t\Big)$.
\paragraph{Example 2.} Let $T=t$ and $H=0$, then
the infinitesimal generator is
\begin{eqnarray*}
v_2=t\,\frac{\partial}{\partial t} + \sum_{i=1}^3 k_i\,t\,
\frac{\partial}{\partial k_i} + \sum_{j=1}^3 n_j(1-t)\,
\frac{\partial}{\partial n_j},
\end{eqnarray*}
Then, the flows of $v_2$ for various values of parameter $s$ are
$$(t,k_i,n_j)\rightarrow\left(t\,e^s,k_i\,e^{t(e^s-1)},n_j\,e^{s-t(e^s-1)}\right).$$
Also, we have the below invariants
\begin{eqnarray*}
I_i=\ln(k_i)-t,\hspace{1cm}
I_{j+3}=\ln\Big(\frac{n_j}{t}\Big)+t,\hspace{1cm} \mbox{for}\:\:
i,j=1,2,3,
\end{eqnarray*}
and the general solution of Eq.~(\ref{eq:2}) as $S(t,{\bf k},{\bf
n}) = \mu\Big(\ln(k_i)-t,\ln\left(\frac{n_j}{t}\right)+t\Big)$
when $\mu$ is an arbitrary function.
\paragraph{Example 3.} For the case which $T=0$ and
$H=k_1\,e^{-t}$, the infinitesimal generator (\ref{eq:20}) changes
to
\begin{eqnarray*}
v_3=  k_1\,\frac{\partial}{\partial k_1} - \sum_{j=1}^3
n_1\,\frac{\partial}{\partial n_j}.
\end{eqnarray*}
The derived group transformations of $v_3$ for parameter $s$ are
$$(t,k_i,n_j)\rightarrow\left(t,k_1\,e^s,k_2,k_3,n_1\,e^{-s},n_2\,
e^{-s}-n_1+n_2,n_3\,e^{-s}-n_1+n_3\right).$$ Thus, the (modified)
invariants are (we suppose that $k_2,k_3\neq 0$, otherwise
$k_2,k_3$ will be two invariants)
\begin{eqnarray*}\begin{array}{lll}
I_1=t, &\hspace{0.5cm} I_2=k_2 ,&\hspace{0.5cm} I_3=k_3,\\
I_4=\displaystyle{\frac{k_1}{n_1}},&\hspace{0.5cm}
I_5=\displaystyle{\ln(k_1)-\frac{n_2}{n_1}},&\hspace{0.5cm}
I_6=\displaystyle{\ln(k_1)-\frac{n_3}{n_1}},
\end{array}\end{eqnarray*}
and for arbitrary function $\mu$, the general solution has the
form
$$S(t,{\bf k},{\bf n})=
\mu\Big(t,k_2,k_3,\frac{k_1}{n_1},\ln(k_1)-\frac{n_2}{n_1},\ln(k_1)-\frac{n_3}{n_1}\Big).$$
\paragraph{Example 4.} If we suppose $T=t$ and
$H=(k_1+k_2+k_3)\,e^{-t}$, then we have the following vector field
\begin{eqnarray*}
v_4=t\frac{\partial}{\partial t} + t\,(2\,k_1+k_2+k_3)
\,\frac{\partial}{\partial k_1}+ t\,\sum_{i=2}^3
k_i\,\frac{\partial}{\partial k_i} + \sum_{j=1}^3 \Big(n_j(1 - t)
- n_1\,t\Big)\,\frac{\partial}{\partial n_j},
\end{eqnarray*}
with group transformations of different parameters $s$, that
transform $(t,k_i,n_j)$ to
\begin{eqnarray*}
&&\hspace{-0.7cm}P(s)=\left(t\,e^s,-(k_2+k_3)\,e^{t\,(e^s-1)}+(k_1+k_2+k_3)\,e^{2\,t\,(e^s-1)},k_p\,e^{t\,(e^s-1)},\right.\\
&&\hspace{3.2cm}\left.n_1\,e^{s-2\,t(e^s-1)},n_1\,e^{s-2\,t(e^s-1)}+(n_q-n_1)\,e^{s-t(e^s-1)}\right),
\end{eqnarray*}
where $p,q=2,3$. Its independent invariants are
\begin{eqnarray*}
&&\hspace{-0.9cm}
I_1=(1-2\,t)\,\ln(k_3)-t\,\ln\left(n_1(1-2\,t)\right),
\hspace{1cm} I_2 =2\,t-\ln(2\,k_1+k_2+k_3),\\
&&\hspace{-0.9cm}
I_3=(1-t)\,\ln(k_2)-t\,\ln\left(n_2(1-t)-n_1\,t\right),\hspace{0.5cm}
I_4 = t
- \ln(k_2)\\
&&\hspace{-0.9cm}
I_5=(1-t)\,\ln(k_3)-t\,\ln\left(n_3(1-t)-n_1\,t\right),\hspace{0.5cm}
I_6 = k_2\,k_3^{-1},
\end{eqnarray*}
and hence the general solution of (\ref{eq:2}) in respect to
infinitesimal operator $v_4$ is an arbitrary function of these
invariants. Indeed, if $u=f(t,k_i,n_j)$ be a solution of
Eq.~(\ref{eq:2}) then also is $u=f(P(s))$ for each $s$.

\section{The Contact Symmetry of the Equation}

In continuation, we change the group action and find symmetry
group and invariants of Eq.~(\ref{eq:2}) up to the contact
transformation groups. According to B\"{a}cklund theorem
\cite{Ol}, if the number of dependent variables be greater than
one (like our problem), then each contact transformation is the
prolongation of a point transformation. But in this section, we
directly earn the structure of infinitesimal generators of contact
transformations\par
We suppose that the general form of a contact transformation be as
following
\begin{eqnarray*}\begin{array}{lll}
\widetilde{t}=\phi(t,k_i,n_j,q_r,p_s),&
\widetilde{k}_l=\chi_l(t,k_i,n_j,q_r,p_s),&
\widetilde{n}_m=\psi_m(t,k_i,n_j,q_r,p_s),\\
\widetilde{q}_n=\eta_n(t,k_i,n_j,q_r,p_s),&
\widetilde{p}_u=\zeta_u(t,k_i,n_j,q_r,p_s),&
\end{array}\end{eqnarray*}
where $i,j,l,m,n$ and $u$ varies between 1 and 6; and
$\phi,\chi_l,\psi_m,\eta_n$ and $\zeta_u$ are arbitrary smooth
functions.
In this case of group action, an infinitesimal generator which is
a vector field in $J^1({\Bbb R},{\Bbb R}^6)$, has the following
general form
\begin{eqnarray*}
v:=T\,\frac{\partial }{\partial t}+\sum_{i=1}^3
\Big[K_i\frac{\partial }{\partial k_i}+ N_i\frac{\partial
}{\partial n_i}+ Q_i\frac{\partial }{\partial q_i}+
P_i\frac{\partial }{\partial p_i} \Big],
\end{eqnarray*}
for arbitrary smooth functions $T, K_l, N_m, Q_m, P_u$ ($l=1,2$
and $1\leq m, n, u\leq3$). \par
Since our computations are done in 1--jet space, so we do not need
to lift $v$ to higher jet spaces and hence we act $v$ (itself) on
the Eq.~(\ref{eq:2}), then we find the following relation
\begin{eqnarray*}
&&\hspace{-0.7cm} \sum_i [n_i(Q_i-K_i)+ N_i(q_i-k_i)]=0.
\end{eqnarray*}
Since ${\bf n}\neq 0$, so without less of generality, we can
suppose that $n_1\neq 0$, then the solution to this equation for
would be
$$K_1 \!=\! Q_1+\frac{1}{n_1}\Big[\sum_{i=2,3}\,n_i(Q_i-K_i)+
\sum_j N_j(q_j-k_j)\Big].$$
Therefore, the infinitesimal generator which we call it as {\it
contact infinitesimal generator} is in the following form
\begin{eqnarray}
v &=& T\,\frac{\partial }{\partial t} + \sum_{i=2}^3
K_i\Big(\frac{\partial }{\partial
k_i}-\frac{n_i}{n_1}\,\frac{\partial }{\partial k_1}\Big)+
\sum_{j=1}^3 \Big[Q_j\Big(\frac{\partial }{\partial
q_i}+\frac{n_j}{n_1}\,\frac{\partial }{\partial k_1}\Big)\nonumber\\
&& + N_j\Big(\frac{\partial }{\partial
n_j}+\frac{1}{n_1}(q_j-k_j)\,\frac{\partial }{\partial k_1}\Big)+
P_j\,\frac{\partial }{\partial p_j}\Big]. \label{eq:24}
\end{eqnarray}
One may divide the latter form to following vector fields, to
consist a basis for Lie algebra ${\frak g}=\langle v \rangle$ of
contact symmetry group $G$
\begin{eqnarray}\begin{array}{rll}
&\hspace{1cm} v_1 = T\,\displaystyle{\frac{\partial }{\partial
t}}, & v_2 = K_2\Big(\displaystyle{\frac{\partial }{\partial k_2}}
-
\frac{n_2}{n_1}\displaystyle{\frac{\partial }{\partial k_1}}\Big), \\[3mm]
&\hspace{1cm} v_3 = K_3\Big(\displaystyle{\frac{\partial
}{\partial k_3}} - \frac{n_3}{n_1}\displaystyle{\frac{\partial
}{\partial k_1}} \Big), & v_4 =
Q_1\Big(\displaystyle{\frac{\partial
}{\partial q_1}} + \frac{\partial}{\partial k_1} \Big), \\[3mm]
&\hspace{1cm} v_5 = Q_2\Big(\displaystyle{\frac{\partial
}{\partial q_2}} + \displaystyle{\frac{n_2}{n_1}\frac{\partial
}{\partial k_1}}\Big), & v_6 =
Q_3\Big(\displaystyle{\frac{\partial
}{\partial q_3}}+\displaystyle{\frac{n_3}{n_1}}\displaystyle{\frac{\partial }{\partial k_1}}\Big), \\[3mm]
&\hspace{1cm} v_7 = N_1\Big(\displaystyle{\frac{\partial
}{\partial n_1}} +
\frac{1}{n_1}(q_1-k_1)\,\displaystyle{\frac{\partial }{\partial
k_1}}\Big), & v_8 = N_2\Big(\displaystyle{\frac{\partial
}{\partial n_2}} + \frac{1}{n_1}(q_2-k_2)\,
\displaystyle{\frac{\partial }{\partial k_1}}\Big),\\[3mm]
&\hspace{1cm} v_9 = N_3\Big(\displaystyle{\frac{\partial
}{\partial n_3}} +
\frac{1}{n_1}(q_3-k_3)\,\displaystyle{\frac{\partial }{\partial
k_1}}\Big),
& v_{10} = P_1\,\displaystyle{\frac{\partial }{\partial p_1}},\\[3mm]
&\hspace{1cm} v_{11} = P_2\,\displaystyle{\frac{\partial
}{\partial p_2}},& v_{12} = P_3\,\displaystyle{\frac{\partial
}{\partial p_3}}.
\end{array}\hspace{-1cm}\label{eq:25}
\end{eqnarray}

The commutators $[v_i,v_j]$ for $1\leq i,j\leq 12$ are linear
combinations of $v_i$ themself, and hence these vector fields
construct a basis for Lie algebra ${\frak g}$ of contact symmetry
group $G$. The commutator table is given in Table 2. In this
table, when commutator of two vector fields has a part in the form
of a $v_i$, then we used $v_i$ instead of it. As is indicated in
this table, for each $1\leq i,j\leq 12$, the Lie bracket of $v_i$
and $v_j$ has two parts, one part in the form of $v_i$, and
another part in the form of $v_j$.
\begin{table}
\begin{eqnarray*}\begin{array}{l|l l }
\hline\hline [\,,\,] &\hspace{1cm} v_i      &\hspace{2cm} v_j
\hspace{1cm} \\ \hline
v_i & \hspace{1cm} 0        &\hspace{2cm} v_i + v_j   \hspace{1cm}\\[1mm]
v_j &\hspace{1cm} -v_i-v_j &\hspace{2cm} 0           \hspace{1cm} \\[1mm]
\hline\hline
\end{array}\end{eqnarray*}
\centering{ \caption{The commutators table provided by contact
symmetry.}}\label{table:2}
\end{table}
\paragraph{Theorem 2.} {\em The contact symmetry group of
(\ref{eq:2}) is an infinite dimensional Lie algebra and its Lie
algebra is generated by contact infinitesimal operators
(\ref{eq:25}) with the commutators table 2.}
%
%%%%%%%%%%%%%%%%%%%%%%%%%%%%%%%%%%%%%%%%%%%%%%%%%%%%%%%%%%%%%%
\section{Conclusion}
A symmetry analysis for a new form of the vortex mode equation led
to find the structure of point and contact infinitesimal
generators as well as fundamental invariants of the new equation.
In addition a form of general solutions implied by these
invariants was obtained. Also we presented some examples for the
point transformation case which tend to a precise determination of
related symmetry groups. In the special case of our problem, the
contact and point symmetry group of the vortex mode equation were
both found to be infinite dimensional Lie groups when the normal
vector to wave front is not necessarily unit.

%%%%%%%%%%%%%%%%%%%%%%%%%%%%%%%%%%%%%%%%%%%%%%%%%%%%%%%%%%%%%%%%%%%%%%%%%%%%%%%%%%%%%%%%%%%%%%%%%%%%%%%%%%%%

\vspace{0.5cm}

\noindent{\it Authors' address:}\\

\noindent Mehdi Nadjafikhah, Ali Mahdipour--Shirayeh\\
School of Mathematics,\\ Iran University of Science and
Technology,\\ Narmak, Tehran 16846 -13114, Iran.\\
E-mail: m\_nadjafikhah@iust.ac.ir, mahdipour@iust.ac.ir
\end{document}